\documentclass[10pt]{amsart}

\usepackage{amsmath}
\usepackage{amsthm}
\usepackage{amsfonts}
\usepackage{amssymb,latexsym}
\usepackage[all]{xy}
\usepackage{graphicx}
\usepackage[mathscr]{eucal}

\theoremstyle{plain}
    \newtheorem{thm}{Theorem}[section]
    \newtheorem{prop}[thm]{Proposition}
    \newtheorem{lemma}[thm]{Lemma}
    
    \newtheorem{cor}[thm]{Corollary}

\theoremstyle{definition}
    \newtheorem{defn}[thm]{Definition}

\theoremstyle{remark}
    \newtheorem{rem}[thm]{Remark}
    \newtheorem{example}[thm]{Example}

\newtheorem*{namedtheorem}{\theoremname}
\newcommand{\theoremname}{testing}

\newcommand{\A}{\mathcal A}
\newcommand{\B}{\mathcal B}
\newcommand{\C}{\mathcal C}
\newcommand{\D}{\mathcal D}
\newcommand{\E}{\mathcal E}
\newcommand{\F}{\mathcal F}

\newcommand{\T}{\mathcal T}

\newcommand{\stk}[1]{\stackrel{#1}{\longrightarrow}}
\newcommand{\la}{\ensuremath{\langle}}
\newcommand{\ra}{\ensuremath{\rangle}}

\newcommand{\ints}{\ensuremath{\mathbb{Z}}}

\newcommand{\ftwo}{\mathbb{F}_2}

\newcommand{\rar}{\ensuremath{\rightarrow}}

\newcommand{\n}{\noindent}

\DeclareMathOperator{\Supp}{Supp}
\DeclareMathOperator{\Spec}{Spec}
\DeclareMathOperator{\Hom}{Hom}

\DeclareMathOperator{\proj}{proj}
\DeclareMathOperator{\Proj}{Proj}
\DeclareMathOperator{\stmod}{stmod}
\DeclareMathOperator{\StMod}{StMod}

\DeclareMathOperator{\Sq}{Sq}
\DeclareMathOperator{\RHom}{RHom}

\DeclareMathOperator{\Ext}{Ext}

\newcommand{\FT}{\mathcal{F}_{tor}}
\DeclareMathOperator{\Cone}{Cone}
\usepackage{verbatim}

\begin{document} 
%
\title{Krull-Schmidt decompositions for thick subcategories}
\date{\today}

\author{Sunil K. Chebolu} 
\address {Department of Mathematics \\
          University of Western Ontario \\
          London, ON. N6A 5B7 }
\email{schebolu@uwo.ca}


\thanks{This research was partially supported by McFarlan Fellowship at the University of Washington.}

\begin{abstract}  Following Krause ~\cite{Kr}, we prove Krull-Schmidt type decomposition theorems for thick 
subcategories of  various triangulated categories including the derived categories of rings, Noetherian stable homotopy categories, stable module categories over 
Hopf algebras, and the stable homotopy category of spectra. In all these categories, it is shown that the thick ideals of small objects decompose uniquely into 
indecomposable thick ideals. We also discuss some consequences of these decomposition results. In particular, it is shown that all these decompositions respect $K$-theory.
\end{abstract}

\maketitle


\section{Introduction}

Using ideas from modular representation theory of finite groups, Krause \cite{Kr} proved a Krull-Schmidt type decomposition theorem for thick subcategories of  
the stable module category of a finite group. More precisely, he showed that the thick ideals of $\stmod(KG)$, where $G$ is a finite group and
$K$ is a field, 
decompose uniquely into indecomposable thick 
ideals. In this paper, we show that such decompositions also exist in various other triangulated categories like the derived categories of commutative
rings and the stable homotopy category of spectra. To this end, we first generalise Krause's 
definition to arbitrary triangulated categories. In our triangulated categories, coproducts will be denoted by $\amalg$ and smash products by $\wedge$.
Throughout this paper, all subcategories will be assumed to be full.

\begin{defn}\cite{Kr} Let $\T$ denote a tensor triangulated category.  A thick subcategory $\A$ of $\T$ is a 
\emph{thick ideal}  if $X \wedge Y$ belongs to $\A$ for all $X \in \A$ and all $Y \in \T$.  If $\A$ is a 
thick (ideal) subcategory of $\T$, a family of thick (ideal) subcategories $(\A_i)_{i \in I}$ is a \emph{decomposition} of $\A$ if 
\begin{enumerate}
\item{the objects of $\A$ are the finite coproducts of objects from the $\A_i$, and}
\item{$\A_i \bigcap \A_j = 0$ for all $i \ne j$.}
\end{enumerate}
A decomposition $(\A_i)_{i \in I}$ of $\A$ is denoted by $\A= \underset{i \in I}{\amalg} \A_i$, and we say  that $\A$ is \emph{indecomposable} 
if $\A \ne 0$  and any decomposition $\A = \C \amalg \D$ implies that $\C = 0$ or $\D = 0$. A decomposition  $\A= \underset{i \in I}{\amalg} \A_i$
into thick subcategories is said to be unique if given another other decomposition  $\A= \underset{i \in I}{\amalg} \B_i$, then 
$\A_i \cong \B_i$ up to a permutation of the indices. 
A \emph{Krull-Schmidt decomposition} of a thick (ideal) subcategory $\A$ is a unique decomposition $\A = \amalg_{i \in I} \A_i$ where all the $\A_i$ are 
indecomposable thick (ideal) subcategories.  
\end{defn}

The above definition reminds one of the classical Krull-Schmidt theorem which states that any finite length module admits a unique direct 
sum decomposition into indecomposable modules. In fact this is the prototypical example on which the results of this papers are based.
Since all thick subcategories studied in this paper consist of compact objects in some triangulated category, 
it is reasonable to call these decompositions as Krull-Schmidt decompositions.

Before we discuss our main results, we should point out the difference between our approach
and the one adopted by Krause. Our approach is very simple and relies heavily on the existing classifications of the lattices of thick subcategories, and hence our results
have a lattice theoretic flavour;  we obtain decompositions for thick subcategories by decomposing their corresponding geometric subsets given by the thick subcategory 
theorems. 
This  approach has the advantage that it applies to a wide range of categories. 
Krause, on the other hand, uses Rickard's idempotent module construction \cite{Ri}  and the theory of endofiniteness \cite[Section 1]{Kr} to decompose the thick
ideals of $\stmod(KG)$; see section \ref{sec:henning} for details.
This approach also has its  advantages. For example, using this approach, Krause arrives at a new
proof (one that does not use Segal's conjecture) of the fact that the classifying space $BG$ of a group $G$ splits stably.

We now discuss the main results in this paper. The first one is about the effect of  $K$-theory of a thick subcategory under a Krull-Schmidt decomposition.

\begin{thm}
Let $\T$ denote a triangulated category and let $\C$ be an essentially small thick subcategory of $\T$. If $\underset{i\in I}{\amalg}{\C_i}$ is a
Krull-Schmidt decomposition of $\C$, then under some mild assumptions (see theorem \ref{th:main}),
there is a natural isomorphism in $K$-theory,
\[ \underset{i \in I}{\bigoplus}\;{K_0(\C_i)} \cong \; K_0(\C).\]
\end{thm}

Our motivation for studying $K$-theory for thick subcategories comes from the problem of classifying the triangulated subcategories in a given
triangulated category. Thomason ~\cite{Th} proved that if $\T$ is
any essentially small triangulated category, then the subgroups of $\T$ are in bijection with the dense triangulated subcategories of $\T$. (A triangulated subcategory
$\A$ is \emph{dense} in $\T$ if $\T$ is the smallest thick subcategory that contains $\A$.) Since every triangulated subcategory $\A$ of $\T$ 
is dense in a unique thick subcategory (namely, the one obtained by taking the intersection of all the thick subcategories of $\T$ that contain $\A$), 
one can classify all the triangulated subcategories of $\T$ by computing the Grothendieck groups of all the thick subcategories of $\T$.  Towards this, 
the above theorem tells us that we can restrict ourselves to indecomposable thick subcategories. For the interested reader, we refer to
\cite{Cheb06a} where we use this brilliant recipe of Thomason to study classifications of triangulated subcategories of finite spectra and perfect complexes.

We use Landsburg's criterion  (lemma ~\ref{le:landsburg}) and some general lemmas about triangulated categories which are 
developed in the Section 2 to prove the above theorem. It applies to all the decompositions we study in this paper.

The next theorem deals with Krull-Schmidt decompositions in the derived categories of rings and the stable homotopy category of spectra.

\begin{thm}
The thick subcategories of finite spectra and those of perfect complexes in the derived category of a Noetherian ring admit Krull-Schmidt decompositions.
Conversely, in both these cases, given any collection of thick subcategories $(\C_i)_{i \in I}$ 
such that $\C_i \bigcap \C_j = 0$ for all $i \ne j$, there exists a unique thick subcategory $\C$ such that $\C= \amalg_{i \in I} \C_i$. 
\end{thm}

As mentioned above, in proving these decomposition theorems, we make good use of the thick subcategory theorems of Hopkins-Smith ~\cite{hs} and Hopkins-Neeman ~\cite{Ne}.
In the derived category of a Noetherian ring $R$, the Hopkins-Neeman result states that the thick subcategories of perfect complexes are in bijection with 
the \emph{specialisation closed} (subset of $\Spec(R)$ that are a union of closed sets). Given a thick subcategory $\T_S$ (corresponding to a specialisation closed
subset $S$ under this bijection),  we define a  graph $G_S$ as follows:  vertices are the minimal primes contained in $S$, and vertices $p$ and $q$ are 
adjacent if $V(p) \cap V(q) \ne \emptyset.$  It is shown that the indecomposable pieces 
that constitute $\T_S$ in a Krull-Schmidt decomposition correspond precisely to the connected components of this graph.
This decomposition theorem gives the following interesting algebraic result; see corollary \ref{cor:siyengar} for a stronger result.

\begin{cor}   Let $X$ be a perfect complex over a Noetherian ring $R$. Then there exists a unique decomposition 
\[ X \cong \underset{i \in I}{\bigoplus}\; X_i\]
such that  the supports of the $X_i$ are pairwise disjoint and indecomposable.
\end{cor}

We also generalise the above results on Krull-Schmidt decomposition to more general
Noetherian stable homotopy categories. We use the language and results from ~\cite{mps} to achieve this generalisation.
In particular, this generalisation gives Krull-Schmidt decompositions in some categories arising in modular representation theory.
We summarise these results in the next theorem. 

\begin{thm} Let $B$ denote a finite dimensional graded co-commutative Hopf algebra satisfying the tensor product property (Examples:
group algebras of finite groups and the finite dimensional sub-Hopf algebras of the mod $2$ Steenrod algebra).
Then,
\begin{enumerate} 
\item Every thick  ideal of small objects in $K(\Proj\,B)$ (the chain homotopy category of graded projective $B$-modules) is indecomposable.
\item Every thick ideal of $\stmod(B)$ admits a Krull-Schmidt decomposition.
\end{enumerate}
\end{thm}

As mentioned above, Krause ~\cite{Kr} proved that thick ideals of $\stmod(KG)$ ($G$ a finite group) admit Krull-Schmidt decompositions.
It is known \cite{bcr2} that the group algebras of finite groups satisfy the tensor product property and therefore part (2) of the above theorem 
is a generalisation of Krause's result.

For the interested reader we mention two recent preprints \cite{bal2} and \cite{buankrausesolberg} as good companions to the present
paper. While the results in these papers do not overlap with ours, nevertheless they address similar issues. As we understand, both these
preprints were motivated by  the notion of  spectrum for a tensor triangulated category which was introduced by Balmer in \cite{bal1}. Using the 
theory developed  in \cite{bal1}, we can easily generalise and state a Krull-Schmidt decomposition theorem for thick ideals of a tensor triangulated category 
with suitable assumptions. However, it is not clear to us if such a generalisation would cover any other example that is not already discussed in this paper. This
is because in all the examples that we know, computing the spectrum of a triangulated category, in the sense of Balmer \cite{bal1}, makes use of a known 
classification of the thick subcategories for the category in question. Therefore we prefer not to get into that generalisation at this point.

The paper is organised as follows. In section  \ref{sec:ktheory}, after developing some necessary lemmas about triangulated categories and thick subcategories, 
we show that Krull-Schmidt decompositions respect $K$-theory under some mild hypothesis. In the later sections
we  study Krull-Schmidt decompositions in various triangulated categories: we review Krause's result for the stable module categories of group algebras in section
\ref{sec:henning}, derived categories of rings in section \ref{sec:KS-derivedcategory}, Noetherian stable homotopy categories in section \ref{sec:nsht},
stable module categories over Hopf algebras in section \ref{sec:KS-hopfalgebras}, and finally the stable homotopy category of spectra in section \ref{sec:KS-finitespectra}.
We end the paper with a few questions which ask for further extensions of these decompositions. \\

\n
\textbf{Acknowledgements:} I would like to thank Srikanth  Iyengar, Henning Krause and John Palmieri for many interesting and fruitful discussions on the subject, 
and the referee for some helpful comments on this paper.

\section{$K$-theory for thick subcategories} \label{sec:ktheory}

In this section we study how the $K$-theory of triangulated categories behaves under a Krull-Schmidt decomposition. 
We begin with some preliminaries on triangulated categories. Let $\T$ denote a  triangulated category that is
essentially small (i.e., a category that has only a set of isomorphism classes of
objects). The \emph{Grothendieck group} $ K_0(T) $ is defined to
be the free abelian group on the isomorphism classes of $\T$ modulo
the Euler relations: $[B]=[A]+[C]$, whenever $A\rightarrow B
\rightarrow C \rightarrow \Sigma A $ is an exact triangle in $\T$
(here $[X]$ denotes the element in the Grothendieck group  that is
represented by the isomorphism class of the object $X$). This is
clearly an Abelian group with $[0]$ as the identity element and
$[\Sigma X]$ as the inverse of $[X]$. We always have the identity
$[A]+[B]=[A \amalg B]$ in the Grothendieck group. Also note that any element of $K_0(\T)$ is of the form $[X]$ for some $X$ in $\T$.

Having defined $K$-theory, we now study how it behaves under a Krull-Schmidt decomposition. We begin with some lemmas that will
be needed in studying $K$-theory. We will start with the following extremely useful Lemma due to Landsburg ~\cite{la}. This is a nice criterion for the 
equality of two classes in the Grothendieck group.

\begin{lemma} \cite{la} \label{le:landsburg}  Let $\T$ be an essentially small triangulated category. If
$X$ and $Y$ are objects in $\T$, then $[X]=[Y]$ in $K_0(\T)$ if
and only if there are objects $A, B$, $C$ in  $\T$ and maps such
that there are exact triangles
\[A \stk{f}  B \amalg X \stk{g} C \stk{h} \Sigma A, \]
\[A \stk{f'} B \amalg Y \stk{g'} C \stk{h'} \Sigma A .\]
\end{lemma}

It is well-known that coproducts of exact triangles are exact in any triangulated category; see \cite[Appendix 2, Proposition 10]{mar}. Under some additional hypothesis,
the next lemma gives a converse to this well-known fact.

\begin{lemma} \label{le:sweet} Let $\T$ denote a triangulated category and let 
\[ A \stk{f} B \rar C \rar \Sigma A  \quad \mbox{and} \quad A' \stk{f'} B' \rar C' \rar \Sigma A' \]
be two sequences of maps in $\T$ such that their sum $ A \amalg A' \stk{f \amalg f'} B\amalg B' \rar C \amalg C' \rar \Sigma (A\oplus A') $
is an exact triangle. If we also have
\begin{eqnarray}
\Hom(\Cone (f),C')  = 0 = \Hom(\Cone (f'),C)  \quad \mbox{and} \\
\Hom(C',\Cone (f))  = 0 = \Hom(C, \Cone(f')), 
\end{eqnarray}
then the given sequences of maps are exact triangles.

\begin{proof} Complete the maps $f$ and $f'$ to exact triangles in $\T$: 
\[ A \stk{f} B \rar \Cone(f) \rar \Sigma A  \quad \mbox{and} \quad A' \stk{f'} B' \rar \Cone(f') \rar \Sigma A'. \]
Since the coproduct of exact triangles is exact, adding these two triangles gives another triangle,  
\[ A \amalg A' \stk{f\oplus f'} B\amalg B' \rar \Cone(f) \amalg \Cone(f') \rar \Sigma (A\oplus A'). \]

We know from the axioms for a triangulated category that there is a fill-in map $H$ in the diagram below.
\[
\xymatrix{ 
A \amalg A' \ar[d]^{=} \ar[r]^{f \amalg f'}&  B \amalg B' \ar[d]^{=} \ar[r] & \Cone(f) \amalg \Cone(f') 
\ar@{..>}[d]^{H} \ar[r] & \Sigma (A \amalg A') \ar[d]^{=} \\ 
A \amalg A' \ar[r]^{f \amalg f'} & B \amalg B' \ar[r] & C \amalg C' \ar[r] &\Sigma(A \amalg A')
 }
\]
Note that three out of the four vertical maps in the above diagram are isomorphisms and therefore so is $H$; see \cite[Appendix 2, Proposition 6]{mar}.
Now the hypothesis 
\[\Hom(\Cone (f), C') = \Hom(\Cone(f'),C)= 0\] 
implies that $H = h \amalg h'$. So we have  $h \amalg h' : \Cone(f) \amalg \Cone(f') \rar C \amalg C'$ is an isomorphism.
The hypothesis 
\[ \Hom(C' , \Cone(f)) = \Hom(C, \Cone(f')) = 0\]
implies that the inverse $G$ of this isomorphism is of the form $ g \amalg g'$.
This forces both $h$ and $h'$ to be isomorphisms and hence $C \cong \Cone(f)$ and $C' \cong \Cone(f')$. 
Since exact triangles in $\T$ form a replete class, the two sequences of maps under consideration are exact triangles.
\end{proof}

\end{lemma}

\begin{lemma} \label{le:short} Let $\C$ and $\D$ be thick ideals in $\T$. If $\C \cap \D = 0$, then $\C \wedge \D = 0$.
($\C \wedge \D$ is the full subcategory of objects of the form $X \wedge Y$ where $X \in \C$ and $Y \in \D$.)
\end{lemma}
\begin{proof} Clear, since $\C \wedge \D \subseteq \C \cap \D = 0$. 
\end{proof}

We are now ready to state and prove our main theorem on $K$-theory for thick subcategories.

\begin{thm} \label{th:main} Let $\A = \underset{i \in I}{\coprod} \A_i$ be a Krull-Schmidt decomposition of an essentially small thick (ideal) subcategory 
$\A$ 
in a triangulated category $\T$. If $\A$ is thick (not necessarily ideal), assume that (1) holds; if $\A$ is a thick ideal, assume any one of the following three 
conditions.
\begin{enumerate} 
\item{$\Hom(\A_i,\A_j)=0$ for all $i \ne j$. }
\item{For all $i$ in $I$, there exists an object $S_i$ in $\A_i$ such that $S_i \wedge X = X$ for all $X \in \A_i$.} 
\item{For all $i$ in $I$, there exists an object $S_i$ in $\T$ such that $S_i \wedge X = X$ for all $X \in \A_i$, and whenever $j \ne i$, $S_i \wedge X = 0$ 
for all $X \in \A_j$.}
\end{enumerate}
Then the inclusion functors  $\A_i \hookrightarrow \A$ give rise to an isomorphism
\[ \Xi: \underset{i \in I}{\bigoplus} \; K_0(\A_i) \cong K_0(\A). \]
\end{thm}

\begin{proof} The inclusion functors  $\A_i \hookrightarrow \A$ induce maps on the Grothendieck groups which can be assembled 
to obtain the map  $\Xi: \underset{i \in I}{\bigoplus} K_0(\A_i) \rar K_0(\A)$ that we want to show is an isomorphism. Showing that 
$\Xi$ is 
surjective is easy: note that every element in  $K_0(\A)$ is of the form $[X]$ for some $X \in \A$. Since the family  
$(\A_i)_{i \in I}$ of thick subcategories is a decomposition for $\A$, we can express $X$ as a finite coproduct; $X = \coprod X_i$ with $X_i \in \A_i$. This 
gives $[X] = [\coprod X_i] = \Sigma_{i \in I} [X_i]$. The last quantity is clearly the image of $\bigoplus [X_i]$ under the map $\Xi$ 
and 
therefore $\Xi$ is surjective. 

To show that $\Xi$ is injective, we use Landsburg's criterion (Lemma \ref{le:landsburg}).  Suppose $\Xi(\bigoplus [X_i]) = 0$. Then, 
as noted above,  $\Xi(\bigoplus [X_i]) = \Sigma [X_i] = [\coprod X_i]$ and hence 
$[\coprod X_i]=[0]$. This now gives, by  Landsburg's criterion, two exact triangles in $\A$: 
\begin{subequations} \label{triangle} 
\begin{gather}
 A \rar B \amalg \left(\underset{j \in I}{\coprod} X_j \right) \rar C \rar \Sigma A  \label{trianglea} \\
A \rar B \rar C \rar \Sigma A  \label{triangleb} 
\end{gather}
\end{subequations}
We want to show that $\bigoplus [X_i] = 0$ or equivalently $[X_i]=0$ for all $i \in I$. 
First assume that condition (1) holds.  Then consider the decompositions $A = \coprod A_i$, $B = \coprod B_i$, and $C = \coprod C_i$
(which exist because $\A = \underset{i \in I}{\coprod} \A_i$ ). Substituting these decompositions
in the above triangles \eqref{trianglea} and \eqref{triangleb} gives
\[\coprod A_i \rar \coprod \left( B_i \amalg X_i \right) \rar \coprod C_i \rar \Sigma \coprod A_i, \]
\[\coprod A_i \rar \coprod B_i \rar \coprod C_i \rar \Sigma \coprod A_i. \]
Now our assumption (1) together with lemma \ref{le:sweet} will enable us to split these two triangles into exact triangles in $\A_i$. 
So for each $i \in I$, we get exact triangles 
\[A_i \rar B_i \amalg X_i \rar C_i \rar \Sigma A_i,\]
\[A_i \rar B_i \rar C_i \rar \Sigma A_i.\]
in $\A_i$. This implies (by Landsburg's criterion) that $[X_i] = 0$ in $K_0(\A_i)$. So $\Xi$ is injective if condition (1) holds.

Now assume that $\A$ is a thick ideal and that for each fixed $i \in I$ either (2) or (3) holds. Smash the above triangles \eqref{trianglea} and \eqref{triangleb} 
with $S_i$ 
to get two exact triangles in $\T$:
\[A \wedge S_i \rar (B \wedge S_i) \amalg \left(\underset{j \in I}{\coprod}( X_j \wedge S_i) \right) \rar C \wedge S_i \rar \Sigma( A 
\wedge S_i)\]
\[A \wedge S_i \rar B \wedge S_i \rar C \wedge S_i  \rar \Sigma (A \wedge S_i)\]
It is easily seen that these triangles are in fact triangles in $\A_i$: This is trivial if condition (2) holds (because $\A_i$ are ideals and 
$S_i \in \A_i$). If condition (3) holds, write
$A = \coprod A_i$ with $A_i \in \A_i$, then $A \wedge S_i = \coprod (A_i \wedge S_i) = A_i \in \A_i$. Similarly $B \wedge S_i$ and $C 
\wedge S_i$ also belong to $\A_i$.

Now we claim that $\underset{j \in I}{\coprod} (X_j \wedge S_i) = X_i$. If (2) holds, then by the lemma \ref{le:short}, we get
$X_j \wedge S_i = 0$ whenever $i \ne j$, and since $S_i$ is a unit for  $\A_i$, $X_i \wedge S_i = X_i$. If (3) holds, this is obviously true.
So in both cases (conditions (2) and (3)), the above triangles can be simplified to obtain the following triangles in $\A_i$.
\[A \wedge S_i \rar (B \wedge S_i) \amalg X_i  \rar C \wedge S_i \rar \Sigma( A \wedge S_i),\]
\[A \wedge S_i \rar B \wedge S_i \rar C \wedge S_i  \rar \Sigma (A \wedge S_i).\]
This implies that $[X_i] = 0$ in $K_0(\A_i)$.  This show that $\Xi$ is injective, completing the proof of the theorem.
\end{proof}

Here is another crucial lemma for studying Krull-Schmidt decompositions.

\begin{lemma} \label{le:unclear} Let $\T$ be a triangulated category and let $\A$ and $\B$ be two thick (ideal) subcategories of $\T$. 
If $\;\Hom(\A, \B) = 0 = \Hom(\B, \A)$, then the full subcategory $\A \amalg \B$, consisting of objects of the form $A \amalg B$ with 
$A \in \A$ and $B \in B$, is a thick (ideal) subcategory of $\T$. 
\end{lemma}

\begin{proof} The key observation here is that every map $H: A \amalg B \rar A' \amalg B'$  in $\A \amalg \B$ is forced by the given hypothesis
to be of the form $f \amalg g$.  It is clear that $\A \amalg \B$ satisfies the ideal condition. To see that $\A \amalg \B$ satisfies
the 2 out of 3 condition, start with a map $H$ as above and complete it to a triangle. Then we have the following diagram where the rows
are triangles in $\T$. (The bottom row is the coproduct of two triangles in $\T$.)
\[
\xymatrix{ 
A \amalg B \ar[d]^{=} \ar[r]^{f \amalg g}&  A' \amalg B' \ar[d]^{=} \ar[r] & \Cone(f \amalg g) 
\ar@{..>}[d]^{H} \ar[r] & \Sigma (A \amalg B) \ar[d]^{=} \\ 
A \amalg B \ar[r]^{f \amalg g} & A' \amalg B' \ar[r] & \Cone(f) \amalg \Cone(g) \ar[r] &\Sigma(A \amalg A')
}
\]
There exits a fill-in map $H$ which turns out to be an isomorphism as before. Therefore $\A \amalg \B$ is a triangulated subcategory.

It remains to show thickness, i.e., $\A \, \amalg \, \B$ is closed under retractions. Consider a retraction map $e: A \amalg B \rar A \amalg B$
(so $e^2=e$). Since $e = a \amalg b$, the equation $e^2 =e$ implies $(a \amalg b)^2 = a^2 \amalg b^2 = a \amalg b$. This shows
both $a$ and $b$ are retractions. So we are done.
\end{proof}

Having developed all the necessary tools, we now turn our attention to Krull-Schmidt decompositions for thick subcategories.
We begin with Krause's decomposition result \cite{Kr} for the stable module category  in the next section.

\section{Stable module categories over group algebras} \label{sec:henning}

Consider the stable module category $\StMod(KG)$, where $G$ is a finite group and $K$ is some field. The objects of this
category are the (right) $KG$-modules and morphisms are equivalence classes of $KG$-module homomorphisms where two homomorphisms are equivalent 
if and only if their difference factors through a projective module. This category is well-known to be a triangulated category and has a well-defined 
tensor product (ordinary tensor product of $K$-vector spaces with the diagonal $G$-action) which makes it into a tensor triangulated category. The full subcategory of 
small objects in $\StMod(KG)$  is equivalent to the full subcategory consisting of finitely generated $KG$-modules and is denoted by $\stmod(KG)$. 
The main theorem of \cite{Kr} then states:

\begin{thm} \cite{Kr}
Every thick ideal $\A$ in $\stmod(KG)$ decomposes uniquely into indecomposable thick ideals $(\A_i)_{i \in I}$; $\A = \coprod_{i \in I} 
\A_i.$ 
Conversely given thick  ideals $(\A_i)_{i \in I}$ such that $\A_i \cap \A_j = 0$ for all $i \ne j \in I$, there exists a thick ideal 
$\A$ such that $\A = \coprod_{i \in I} \A_i$.  
\end{thm}

We will briefly outline how Krause arrives at this decomposition. The key idea is to consider the Bousfield localisation with respect to the localising 
subcategory  $\A^{\amalg}$ generated by $\A$. The inclusion $\A^{\amalg} \hookrightarrow \StMod(KG)$ has a right adjoint 
$e: \StMod(KG) \rar {\A}^{\amalg}$. In fact, $e(M)$ is just the fibre of the Bousfield localisation map $M \rar M_{A^{\amalg}} $. Thus for each $KG$-module 
$M$, there is a natural triangle in $\StMod(KG)$:
\[e(M) \stk{\epsilon} M \rar M_{\A^{\amalg}} \rar \Sigma e(M).\]
In particular, when $M$ is the trivial representation $K$, we get a module $e(K)$, which we denote by $E_\A$.
The $KG$-module $E_\A$ associated to the thick ideal $\A$ in this way is an idempotent module ($E_\A \otimes E_\A = E_\A$). (These modules were introduced by
Rickard and they proved to be very useful objects in modular representation theory.) Krause then shows that this 
idempotent module is an endofinite object (see \cite[Definition 1.1]{Kr}) in
$\A^{\amalg}$ and hence admits a splitting into indecomposable modules: $E_\A = \bigoplus_{i \in I} E_i$ with $E_i \otimes E_j =0$ for 
$i \ne j$. 
Let $(\epsilon_i): \coprod_i E_i \rar K$ be the decomposition of $\epsilon: E_\A \rar K$ and define $\A_i$ to be the collection of all 
modules $X$ in 
$\stmod(KG)$ such that $\epsilon_i \otimes X$ is an isomorphism. It then follows that $\A = \coprod_i \A_i$ is the Krull-Schmidt 
decomposition for $\A$.
See \cite{Kr} for more details.

We draw the following corollary by applying theorem \ref{th:main} to the above decomposition.

\begin{prop} \label{prop:slick} Let $\A$ be a thick ideal of $\stmod(KG)$ and let $\amalg \, \A_i$ be the Krull-Schmidt decomposition of $\A$. Then,
\[ K_0(A) \cong \underset{i \in I}{\bigoplus} \;K_0(\A_i).\]
\end{prop}
\begin{proof} We show that condition (1) of theorem \ref{th:main} is satisfied. This is done in \cite[Lemma 2.4]{Kr} but we include it here 
for the reader's convenience. For any two thick ideals  $\A_1$ and $\A_2$ that occur in the decomposition for $\A$, we want to show that 
$\Hom(\A_1, \A_2)=0$. Towards this, consider objects $X \in \A_1$ and $Y \in \A_2$ and note that every map $X \rar Y$ in $\stmod(KG)$ 
factors through $X \otimes \Hom_K(X,Y)$ in the obvious way. So we will be done if we can show that $\Hom_K(X,Y)$ is a projective $KG$-module 
or equivalently a trivial object  when viewed in  $\stmod(KG)$. To see this, first note that 
$\Hom_K(X,Y) \cong X^* \otimes Y$, where $X^*$ denotes the $K$-dual of $X$. But $X^*$ is a retract of $X^* \otimes X \otimes X^*$
(see \cite[Lemma A.2.6]{mps}) and therefore belongs to $\A_1$. This implies that $X^* \otimes Y \in \A_1 \otimes \A_2$. 
The latter is zero by lemma \ref{le:short} and therefore $\Hom_K(X,Y) = X^* \otimes Y$ is zero in $\stmod(KG)$. 
\end{proof}

\section{Derived categories of rings} \label{sec:KS-derivedcategory}

In this section we will prove a Krull-Schmidt theorem for the thick subcategories of perfect complexes over a Noetherian ring.
We work in some level of formality here so that we can generalise our results easily to Noetherian stable homotopy categories in the sense
of Hovey-Palmieri-Strickland \cite{mps}. We begin with some preliminaries.

Recall that a subset $S$ of $\Spec(R)$ is a \emph{thick support} if it is a union of (Zariski) closed sets $S_{\alpha}$ such that $\Spec(R) - S_{\alpha}$ 
is a quasi-compact set. It is an exercise \cite[Page 12, Exercise.17(vii)]{am} to show that $\Spec(R) - S_{\alpha}$ is quasi-compact if and only if 
$S_{\alpha} = V(I_{\alpha})$ for some finitely generated ideal $I_{\alpha}$. ($V(I)$ denotes the collection of all primes that contain $I$.)
Since all ideals in a Noetherian ring are finitely generated, it follows that the thick
supports for Noetherian rings are precisely subsets of $\Spec(R)$ that are a union of closed sets (also known as \emph{specialisation closed subsets}).
It is a well-known fact that the compact objects in the derived category of a ring are the ones that are quasi-isomorphic to perfect complexes (bounded chain complexes of
finitely generated projective modules); see, for instance, \cite[Proposition 9.6]{ch}.
The full subcategory of perfect complexes will be denoted by $D^b(\proj\,R)$. If $X$ is a complex in this category we define
its support as 
\[\Supp(X) = \{ p \in \Spec(R) \,| \, X \otimes_R R_p \ne 0\}.  \]
We now state the celebrated thick subcategory theorem for the derived category.

\begin{thm}\cite{Th} The lattice of thick subcategories in $D^b(\proj\,R)$ is isomorphic to the lattice of thick supports of $\Spec(R).$
Under this isomorphism, a thick support $S$ corresponds to the thick subcategory $\T_S$ consisting of all complexes $X$ such that $\Supp(X) \subseteq S$, 
and a thick subcategory $\C$ corresponds to the thick support $\bigcup_{X \in \C} \Supp(X).$
\end{thm}

Using this theorem, we now work our way to the Krull-Schmidt theorem for thick subcategories of perfect complexes. We begin with some lemmas.

\begin{lemma} \label{le:support} For $X$ in $D^b(\proj\,R)$, let $DX = \RHom(X,R)$ denote its Spanier-Whitehead dual. Then $\Supp(X) = \Supp(DX).$ 
\end{lemma}
\begin{proof} $X$ is a retract of $X \otimes DX \otimes X$ \cite[Lemma A.2.6]{mps}, and therefore $DX$ is a retract of $DX \otimes X \otimes DX$. 
Now it is clear that both $X$ and $DX$ have the same support.
\end{proof}

Following \cite{mps}, we will denote $\Hom_{D^b(\proj\,R)}(\Sigma^* X ,Y)$ by $[X,Y]_*$ and the internal function object 
$\RHom(X,Y)$ by $F(X,Y)$. With these notations, we have the following natural isomorphism \cite{mps}
\[[X \otimes A , B]_* \cong [X , F(A,B)]_* .\] 
This isomorphism gives the following useful lemma.

\begin{lemma} \label{le:cond1} If $X$ and $Y$ are perfect complexes such  that $\Supp(X)$ and $\Supp(Y)$ are disjoint, 
then $[X,Y]_* = 0.$ In particular if $A$ and $B$ are any two disjoint thick supports of $\Spec(R)$, then  $[\T_A, \T_B]_* = 0$.
\end{lemma}

\begin{proof} By Spanier-Whitehead duality, $[X,Y]_* \cong [R \otimes X, Y_*] \cong [R, F(X,Y)]_*$. Therefore $[X,Y]_* = 0 $ if and only if
$F(X,Y) = 0$. But since $X$ is a small object in $D(R)$, $F(X,Y) = DX \otimes Y$ \cite[Appendix A.2]{mps}. So we have to show that 
$DX \otimes Y = 0$. Since $\Supp(DX) = \Supp(X)$  is disjoint with $\Supp(Y)$, given any prime $p$, either 
$p$ is not in $\Supp(DX)$ or it is not in $\Supp(Y)$. In the former case, $R_p \otimes DX = 0$, and in the latter, $R_p \otimes Y =0$.
In either case, we get $DX \otimes Y \otimes R_{p} = 0$. Since $p$ is an arbitrary prime, we get $DX \otimes Y = 0$. 
\end{proof}

By a \emph{thick decomposition} of a thick support $S$, we will mean a decomposition $S = \bigcup S_i$ into thick supports, 
where $S_i \cap S_j = \emptyset$ if $i \ne j$. A thick decomposition is \emph{Krull-Schmidt} if the $S_i$ are nonempty, do not admit
nontrivial thick decompositions, and if any such decomposition is unique up to a permutation of the subsets  $S_i$. We illustrate this with 
a couple of examples.

\begin{example} All rings considered are commutative.
\begin{enumerate}
\item{If $R$ is a PID, then every non-zero prime ideal is a maximal ideal and so the thick supports are:
$\Spec(R)$, and all subsets of the set of maximal ideals. In particular, $\Spec(R)$ is an indecomposable thick support.}
\item{If $R$ is an Artinian ring, it has only finitely many prime ideals and every prime ideal is also a maximal ideal. So every 
non-empty subset $S$ of $\Spec(R)$ is a thick support. It is then clear that $\bigcup_{p \in S} \{ p \}$ is the Krull-Schmidt
decomposition for $S$.}
\end{enumerate}
\end{example}

\n
\emph{Note:} The components of the Krull-Schmidt decomposition of a thick support $S$ are not (in general) the connected 
components of $S$. For example,  $\Spec(\ints) - \{ 0 \}$ is a connected subset of $\Spec(\ints)$; however, it is not
indecomposable as a thick support. In fact, 
\[ \Spec(\ints) - \{ 0 \} = \underset{p \ \text{prime}}{\bigcup}\; \{ (p) \}\]
is its Krull-Schmidt decomposition.

Now we establish the strong connection between the decompositions of thick supports and the decompositions of thick subcategories.

\begin{prop} \label{prop:1} Let $C = A \cup B$ be a thick decomposition of a thick support. Then this induces a decomposition of the associated thick
subcategories:
\[\T_C = \T_A \amalg \T_B.\]
\end{prop}
\begin{proof} Clearly $\T_A \cap \T_B = 0$ because $A$ and $B$ are disjoint by definition. We have to show that the objects of $\T_C$
are coproducts of objects in $\T_A$ and $\T_B$. We give an indirect proof of this statement using the thick subcategory theorem.
Lemma \ref{le:cond1} and lemma \ref{le:unclear} tell us that $\T_A \amalg \T_B$ is a
thick subcategory. So we will be done if we can show that the thick support corresponding to  $\T_A \amalg \T_B$ is $C$.
The thick support corresponding to $\T_A \amalg \T_B$  is  given by $\underset{X \in \T_A \amalg \T_B}{\bigcup} \Supp(X)$.
This clearly contains both $A$ and $B$ and hence their union ($ = C) $.  To see the other inclusion,
just note that $\Supp(a \oplus b) = \Supp(a) \cup \Supp(b) \subseteq A \cup B = C.$
\end{proof}

Now we show that decompositions of thick subcategories gives rise to thick decompositions of the corresponding thick supports.

\begin{prop} \label{prop:2} Let $\C = \A \, \amalg \,\B$ be a decomposition of a thick subcategory $\C$ of $D^b(\proj\,R)$ and let $A$, $B$, and $C$ denote the
corresponding thick supports of these thick subcategories. Then $C = A \cup B$ is a thick decomposition of $C$.
\end{prop} 
\begin{proof} We first observe that the intersection of thick supports is again a thick support: Let $S$ and $T$ be two thick supports of $\Spec(R)$.
Then $S = \bigcup_{\alpha} V(I_{\alpha})$ and $T = \bigcup_{\beta} V(J_{\beta})$ for some finitely  generated ideals $I_{\alpha}$ and $J_{\beta}$
in $R$. Now,
 \[ S \cap T = \bigcup_{\alpha , \beta} \left(V(I_{\alpha}) \cap V(J_{\beta})\right) = \bigcup_{\alpha , \beta} V(I_{\alpha} + J_{\beta}).\] 
Since the sum of two finitely generated ideals is  finitely generated, we conclude that  $S \cap T$ is a thick support.

We now argue that $A$ and $B$ are disjoint. If $A$ and $B$ are not disjoint, then their intersection being a nonempty 
thick support corresponds to a non-zero thick subcategory that is contained in both $\A$ and $\B$. This contradicts the fact that $\A \cap \B = 0$, so $A$ and $B$ 
have to be disjoint. 
It remains to show that $A \cup B$ is $C$.  For this, note that every complex $c \in \C$ splits as $c = a \amalg b$ and recall that 
$\Supp(c) = \Supp(a) \cup \Supp(b)$. It follows that $C = A \cup B.$ 
\end{proof}

Combining proposition \ref{prop:1} and \ref{prop:2} we get the following decomposition result.

\begin{thm}\label{thm:aux} Let $R$ be a commutative ring and let $\T_S$ be the thick subcategory of $D^b(\proj \,R)$ corresponding to a thick support $S$.
Then $\T_S$ admits a Krull-Schmidt decomposition if and only if $S$ admits one.
\end{thm}

Applying our main theorem \ref{th:main} to the above decomposition gives,

\begin{cor} \label{cor:Ktheory} Let  $\T_S  = \coprod_{i \in I} \T_{S_i}$ be a decomposition corresponding to a thick decomposition
$S = \bigcup_{i \in I} S_i$. Then,
\[ K_0(\T_S) \cong \underset{i \in I}{\bigoplus} \; K_0(\T_{S_i}).\]
\end{cor}
\begin{proof} Condition (1) of our main theorem holds here by lemma \ref{le:cond1}.
\end{proof}

The question that remains to be addressed is the following. When do thick supports of $\Spec(R)$ admit 
Krull-Schmidt decompositions? We show that this is always possible if $R$ is Noetherian. 

\begin{prop} \label{prop:graph} Let $R$ be a Noetherian ring and let $S$ be a thick support of $\Spec(R)$. 
Then there exists a Krull-Schmidt decomposition $\bigcup S_i$ for $S$.
\end{prop}

\begin{proof} 
It is well-known that the set of prime ideals in a Noetherian ring satisfies the descending
chain condition \cite[Corollary 11.12]{am}. To start, let $S$ be a thick support in $\Spec(R)$ and let $(p_i)_{i \in I}$ be the 
collection of all minimal elements in $S$ -- i.e., primes $p \in S$ which do not contain any other prime in $S$.  It is now clear 
(using the above fact about Noetherian rings) that every prime $p \in S$ contains a minimal element $p_i \in S$, therefore 
$S = \bigcup_{i \in I} V(p_i)$. (Also note that each  $V(p_i)$ is a closed subset of $\Spec(R)$ and hence a thick support.)  
Now define a graph $G_S$ of $S$ as follows: The  vertices are the minimal primes $(p_i)_{i \in I}$ in $S$, and two 
vertices $p_i$ and $p_j$ are adjacent if and only if $V(p_i) \cap V(p_j) \ne \emptyset$. 
Let $(C_k)_{k \in K}$ be the connected components of this graph and for each $C_k$ define a thick support
\[S_k:= \underset{p_i \in C_k}{\bigcup}{ V(p_i)}.\]
By construction it is clear that $\bigcup S_k$ is a thick decomposition of $S$. It is not hard to see that each $S_k$ is indecomposable.
Finally the uniqueness part: let $\bigcup T_k$ be another Krull-Schmidt decomposition of $S$. 
It can be easily verified that the minimal primes in $T_k$ are precisely the minimal primes of $S$ that are contained in $T_k$. Thus the 
Krull-Schmidt decomposition $\bigcup T_k$ gives a partition of the set of minimal primes in $S$. This partition induces a decomposition of the
thick graph of $S$ into its connected (since each $T_k$ is indecomposable) components. Since the decomposition of a graph into its connected 
components is unique, the uniqueness of Krull-Schmidt decomposition follows.
\end{proof}

\begin{rem} The careful reader will perhaps note that the proof of proposition \ref{prop:graph} makes use of \emph{only} the following two 
conditions on $R$.
\begin{enumerate}
\item{Every open subset of $\Spec(R)$ is compact.}
\item{$\Spec(R)$ satisfies the descending chain condition.}
\end{enumerate}
Therefore the proof generalises to any $R$ which satisfies these two properties. As mentioned above, it is well-known  
that Noetherian rings satisfy these two properties. Moreover, any ring which has finitely
many prime ideals automatically satisfies these properties. Note that there are non-Noetherian rings which have finitely many primes. For example,
the ring \[ R = \mathbb{F}_2 [x_2,x_3,x_4, \cdots]/(x_2^2,x_3^3,x_4^4,\cdots)\] 
is a non-Noetherian ring with only one prime ideal. So everything that we are going to state for the remainder of this section will work for rings
which have these two properties, but we state our results only for Noetherian rings for simplicity and cognitive reasons.
\end{rem}

After all that work, the following theorem is now obvious.

\begin{thm} \label{th:mainalg} If $R$ is a Noetherian ring, then every thick subcategory of $D^b(\proj\,R)$ admits a unique Krull-Schmidt decomposition.
Conversely, given any collection of thick subcategories $(\C_i)_{i \in I}$ in $D^b(\proj\,R)$
such that $\C_i \cap \C_j = 0$ for all $i \ne j$, there exists a thick subcategory $\C$ such that $\C= \coprod_{i \in I} \C_i$.
\end{thm}
\begin{proof} The first part follows from theorem \ref{thm:aux} and proposition \ref{prop:graph}. For the second part, define $\C$ be the
full subcategory of all finite coproducts of objects from the $\C_i$. Thickness of $\C$ follows by combining  theorem \ref{thm:aux}, lemma \ref{le:cond1}, 
and lemma \ref{le:unclear}.
\end{proof}

The next corollary is a categorical characterisation of local rings among Noetherian rings.

\begin{cor} A Noetherian ring $R$ is local if and only if every thick subcategory of $D^b(\proj\,R)$ is indecomposable.
\end{cor}
\begin{proof} If $R$ is local, then it is clear that every thick support of $\Spec(R)$ contains the unique maximal ideal. In particular,
we cannot have two non-empty disjoint thick supports and therefore every thick subcategory of $D^b(\proj\,R)$ is indecomposable.
Conversely, if $R$ is not local, it is clear that the thick subcategory supported on
the set of maximal ideals is decomposable.
\end{proof}

Theorem \ref{th:mainalg} also gives the following interesting algebraic results.

\begin{cor} \label{cor:purealg}  Let $X$ be a perfect complex over a Noetherian ring $R$. Then $X$ admits a unique splitting into perfect complexes,
\[ X \cong \underset{i \in I}{\bigoplus} X_i\]
such that the supports of the $X_i$ are pairwise disjoint and indecomposable.
\end{cor}

\begin{proof} Let $\bigcup S_i$ be the Krull-Schmidt decomposition of $\Supp(X)$ (existence is guaranteed by proposition \ref{prop:graph}). 
It is clear from  theorem \ref{th:mainalg} that $X$ admits a splitting ($X \cong \oplus X_i$) where $\Supp(X_i) \subseteq S_i$. To see that we have equality  
in this inclusion, observe that 
\[ \Supp (X) = \Supp \left( \underset{i \in I}{\bigoplus} X_i \right) = \underset{i \in I}{\bigcup} \Supp(X_i) 
\subseteq \underset{i \in I}{\bigcup} S_i = \Supp(X) .\]
Uniqueness: If $X$ admits another decomposition $\oplus Y_i$ as above, then by proposition \ref{prop:graph} we know that both $\bigcup\, \Supp(X_i)$ 
and $\bigcup \, \Supp(Y_i)$ are the same decompositions of $\Supp(X)$.  Lemma  \ref{le:cond1} then implies that $X_i \cong Y_i$ up to a permutation
of the indices. So the decomposition is unique.
\end{proof}

Srikanth Iyengar has pointed out  that this splitting holds in a much more generality. In fact it holds for all complexes in the derived category which have 
bounded and finite homology (i.e. complexes whose homology groups $H_i(-)$ are finitely generated and are zero for all but finitely many $i$.) The full 
subcategory of such complexes will be denoted by $D^f(R)$. 

\begin{cor} \label{cor:siyengar} Let $R$ be a Noetherian ring. Every complex $X$ in $D^f(R)$ admits a splitting 
\[ X \cong \underset{i \in I}{\bigoplus} X_i\]
such that the supports of the $X_i$ are pairwise disjoint and indecomposable.
\end{cor}

\begin{proof} Since $X$ has bounded and finite homology, $\Supp(X)$ is a closed set and therefore proposition \ref{prop:graph} applies and 
we get a decomposition 
\[ \Supp(X) =  \underset{i \in I}{\bigcup} S_i.\]
We now use a result of Neeman \cite{Ne} which states that the lattice of localising subcategories (thick subcategories that are closed under
arbitrary coproducts) is isomorphic to the lattice of all subsets of $\Spec(R)$.
If $\mathcal{L}$ denotes the localising subcategory generated by $X$, then we get a Krull-Schmidt decomposition 
$ \mathcal{L} \cong \coprod \mathcal{L}_i$ that corresponds to the above decomposition of the support of $X$. (This follows exactly as 
in 
the case of thick subcategories.)  The splitting of $X$ as stated in the corollary is now clear. 
\end{proof}

\begin{example} We illustrate the last two corollaries with some examples.
\begin{enumerate}
\item{For every integer $n > 1$, let $M(n)$ denote the Moore complex $0 \rar \mathbb{Z} \stk{n} \mathbb{Z} \rar 0$
and let $p_1^{i_1}p_2^{i_2}\cdots p_k^{i_k}$ be the unique prime factorisation of $n$. It is easy to see that 
$\bigcup_{t=1}^k \{ (p_t)\} $ is the Krull-Schmidt decomposition of $\Supp (M(n))$.
Then it follows that
\[M(n) \cong \bigoplus_{t=1}^k\, M(p_t^{i_t})\] 
is the splitting of $M(n)$  corresponding to this Krull-Schmidt decomposition.}
\item{ Let $R$ be a self-injective Noetherian ring (hence also an Artinian ring). Then it follows (because every prime ideal is also maximal) 
that $\bigcup_{p \in \Spec(R)} \{ p\}$ 
is the Krull-Schmidt decomposition  for $\Supp(R)$.  It can be shown that 
\[ R \cong \bigoplus_{p \in \Spec(R)} E(R/p)\] is the splitting of $R$ corresponding to this Krull-Schmidt decomposition.
(Here $E(R/p)$ denotes the injective hull of $R/p$.)}
\end{enumerate}
\end{example}

\section{Noetherian stable homotopy categories.} \label{sec:nsht}

Motivated by the work in the previous section, we now state a Krull-Schmidt theorem for Noetherian stable homotopy categories.
We use the language and results of \cite{mps} freely. We explain how the proofs of the previous section generalise to give us a
 Krull-Schmidt theorem in this more general setting by invoking the appropriate results from \cite{mps}.
We begin with some definitions and preliminaries from Axiomatic Stable Homotopy Theory \cite{mps}.

\begin{defn} \cite{hp} A \emph{unital algebraic stable homotopy category} is a tensor triangulated category $\mathscr{C}$ with the following 
properties.
\begin{enumerate}
\item{Arbitrary products and coproducts of objects in $\mathscr{C}$ exist.}
\item{$\mathscr{C}$ has a finite set $\mathscr{G}$ of weak generators, i.e., $X \cong 0$ if and only if $[A,X]_*=0$ for all $A \in \mathscr{G}$.}
\item{The unit object $S$ and the objects of $\mathscr{G}$ are small.}
\end{enumerate}
$\mathscr{C}$ is a \emph{Noetherian stable homotopy category} if in addition the following conditions are satisfied,
\begin{itemize}
\item{ $\pi_*(S) := [S,S]_{*}$ is commutative and Noetherian as a bigraded ring. }
\item {For small objects $Y$ and $Z$ of $\C$, $[Y,Z]_{*}$ is a finitely generated module over $[S,S]_{*}$. }
\end{itemize}

 The stable homotopy category of spectra is unital and algebraic but not Noetherian. The derived
category  $D(R)$ of a commutative ring $R$ is a Noetherian stable homotopy category if and only if $R$ is Noetherian
(since $[R,R]_*=R$). $\Spec(\pi_*S)$ will stand for the collection of all homogeneous prime ideals of $\pi_*S$ with the
Zariski topology. 
\end{defn}

Henceforth $\mathscr{C}$  will denote a  bigraded Noetherian stable homotopy category. For each  thick ideal $\mathcal{A}$
of finite objects in $\mathscr{C}$, there is a \emph{finite localisation functor} $L_{\A}$ (also denoted $L_{\A}^f$) on $\mathscr{C}$ whose finite acyclics 
are precisely
the objects of $\A$; see \cite[Theorem 2.3]{hp}. For each bihomogeneous prime ideal $\mathbf{p}$ in $\pi_*(S)$, there are finite localisation
functors $L_{\mathbf{p}}$ and $L_{< \mathbf{p}}$ whose finite acyclics are 
$\{X \mbox{\;finite\;}\;|\; X_{\mathbf{q}} = 0 \;\forall \; \mathbf{q} \subseteq \mathbf{p} \}$ and 
$\{X \mbox{\;finite\;}\;| \; X_{\mathbf{q}} = 0 \;\forall \; \mathbf{q} \subsetneq \mathbf{p}  \}$ respectively. 
These functors define, for each $X \in \mathscr{C}$, a natural exact triangle
\[M_{\mathbf{p}} X \rar L_{\mathbf{p}} X \rar L_{<\mathbf{p}} X \rar \Sigma M_{\mathbf{p}} X.\]
We set $M_{\mathbf{p}} : = M_{\mathbf{p}} S$. See \cite{hp} for more details.

We say that $\mathscr{C}$ has the \emph{tensor product property} if for all bihomogeneous prime ideals $\mathbf{p}$ in $\pi_*(S)$ and all objects 
$X, Y \in \mathscr{C}$,  $M_{\mathbf{p}} \wedge X \wedge Y = 0$ if and only if either $M_{\mathbf{p}} \wedge X = 0$ or  $M_{\mathbf{p}} \wedge Y =0$.

We now mimic the setup for the derived category of a ring. For $X$ in $\mathscr{C}$, define its \emph{support variety} as
\[\Supp(X) = \{\mathbf{p}\;|\; M_{\mathbf{p}} \wedge X \ne 0 \} \subseteq \Spec (\pi_*S),\]
and if $\D$ is any thick subcategory of $\mathscr{C}$, define 
\[\Supp(\D) = \underset{X \in \D}{\bigcup} \Supp(X).\]  
It is a fact \cite[Theorem 6.1.7]{mps} that for $X$ a small object in $\mathscr{C}$,  $\Supp(X)$ is a Zariski closed subset and hence also a thick support of
$\Spec \pi_*(S)$.

We are now ready to state the thick subcategory theorem for $\mathscr{C}$. 

\begin{thm} \label{th:tshps} \cite{hp} Let $\mathscr{C}$ be a Noetherian stable homotopy category satisfying the tensor product property.
Then the lattice of thick ideals of small objects in $\mathscr{C}$ is isomorphic to the the lattice of thick
supports of $\Spec(\pi_*S)$. Under this isomorphism, a thick ideal $\A$ corresponds to $\Supp(\A)$, and
a thick support $T$ of $\Spec(\pi_*S)$ corresponds to  the thick ideal $\{ Z \in \mathscr{C}\;| \;\Supp(Z) \subseteq T\}.$
\end{thm}

Now using this thick subcategory theorem, we can run our proofs from the previous section which generalise to prove theorem \ref{th:grand} below. 
We will  mention the minor changes that are necessary for this generalisation. We begin with an algebraic lemma.

\begin{lemma} \label{le:benson} Let $R$ be a graded commutative Noetherian $K$-algebra. Then $R$ satisfies the descending chain condition on the set 
$\Spec(R)$ of homogeneous prime ideals of $R$.
\end{lemma}

\begin{proof} First note that if char($K$) = 2, then $R$ is a strictly commutative graded ring and the result for such $R$ is well-known; see \cite[Corollary 11.12]{am}.
So assume that $\text{char} (K) \ne 2$. Since every homogeneous prime ideal contains the nilradical, the ring map
\[R \rar R/{\text{nilrad}\, R} \]
induces an order-preserving bijection between $\Spec(R)$ and $\Spec(R/{\text{nilrad}\, R})$. So in view of this remark, we can assume without 
loss of generality that $\text{nilrad} \,R = 0$. But if  $\text{nilrad}\, R = 0$, then $R$ is concentrated in even degrees: If $x$ is an odd 
degree element,  then $x^2 = -x^2$.  This implies $x^2 = 0$ (since char$(K) \ne 2$), or equivalently $x \in \text{nilrad}\, R = 0$.  
This shows that $R$ is concentrated 
in even degrees. In particular it is a strictly commutative graded ring, so we are done.
\end{proof}

\begin{thm} \label{th:grand} Let $\mathscr{C}$ be a Noetherian stable homotopy category satisfying the tensor product property.
Then every thick ideal $\A$ of small objects in $\mathscr{C}$ admits a  Krull-Schmidt decomposition: 
$\A = \coprod_{i \in I} \A_i$. Consequently, 
\[ K_0(\A) \cong \underset{i \in I}{\bigoplus}\;\; K_0(\A_i).\] 
Conversely, if $\A_i$ are thick ideals  of small objects of $\mathscr{C}$ such that $\A_i \cap \A_j = 0$ for all $i \ne j$, then there exists a thick
ideal $\A$ such that $\A = \coprod_{i \in I} \A_i$.  
\end{thm}

\begin{proof} We will explain why all the lemmas and propositions leading up to theorem \ref{th:mainalg} hold in this generality. 
We proceed in order starting from lemma \ref{le:support}.

Lemma \ref{le:support}: This holds for any strongly dualisable object in a closed symmetric monoidal category \cite[Lemma A.2.6]{mps}. 
By \cite[theorem 2.1.3 (d)]{mps}, every small object is strongly dualisable in a unital algebraic stable homotopy category and therefore lemma \ref{le:support} 
holds for small objects of $\mathscr{C}$.

Lemma \ref{le:cond1}:  Proceeding in the same way, it boils down to showing that $DX \wedge Y$ is zero. Since  $DX$ and $Y$ have disjoint 
supports, for each prime $\mathbf{p}$ in $\Spec(\pi_* S)$, either $Y \wedge M_{\mathbf{p}} = 0$ or 
$DX \wedge M_{\mathbf{p}} = 0$. This implies that $DX \wedge Y \wedge M_{\mathbf{p}} = 0$ for all $\mathbf{p}$. By theorem  \cite[Theorem 6.1.9]{mps}, 
we have the following equality of Bousfield classes:
\[ \la S \ra = \underset{\mathbf{p} \,\in \,\Spec(\pi_*S)}{\coprod} \la M_{\mathbf{p}}\ra . \]
This implies $DX \wedge Y = 0$.

Proposition \ref{prop:1}: This goes through verbatim without any changes.

Proposition \ref{prop:2}: This is actually much easier. Thick supports of $\Spec(\pi_*S)$ are just unions of Zariski-closed sets
and therefore it is obvious that intersection of thick supports is thick. The rest of the proof of this proposition follows exactly the same
way. 

Using these lemmas and proposition, theorem \ref{thm:aux} and corollary \ref{cor:Ktheory} follow. Proposition \ref{prop:graph} also holds in the graded 
Noetherian case; see lemma \ref{le:benson}. Thus the first part of theorem follows by combining all these lemmas and propositions.

The converse follows as before by combining theorem \ref{thm:aux}, lemma \ref{le:cond1}, and lemma \ref{le:unclear}.
\end{proof}

\begin{cor}  Let $\mathscr{C}$ be a Noetherian stable homotopy category satisfying the tensor product property. Then every small object $X$ in $\mathscr{C}$ admits a
unique  decomposition
\[ X \cong \underset{i \in I}{\coprod} \; X_i,\]
such that the support varieties of the $X_i$ are pairwise disjoint and indecomposable.
\end{cor}

\begin{rem} The derived category $D(R)$ of a Noetherian ring $R$ is a Noetherian stable homotopy category which satisfies the tensor product property.
(Note that $\pi_*(S)= [R,R]_* = R$.)
So theorem \ref{th:tshps} recovers our Krull-Schmidt theorem for derived categories of Noetherian rings.
\end{rem}

We now give another example of a Noetherian stable homotopy category for 
which theorem \ref{th:grand} applies.

\begin{example} \cite{hp, hp2} 
Let $B$ denote a finite dimensional graded co-commutative Hopf algebra over a field $k$ ($\text{char}(k) > 0$).  Then the category $K(\Proj\,B)$, whose objects are 
unbounded chain complexes of projective $B$-modules and whose morphisms are chain homotopy classes of chain maps, is a Noetherian stable homotopy category.
In this category the unit object is a projective resolution of $k$ by $B$-modules and therefore $\pi_*(S)$ is isomorphic
to $\Ext_B^{*,*}(k,k)$. It is well-known that this ext algebra is a bigraded Noetherian algebra; see \cite{fs} and \cite{wil}.
So theorem \ref{th:grand} applies and we get the following result.
\end{example}

\begin{thm} \label{thm:KSforstable(B)} Let $B$ denote a finite dimensional graded co-commutative Hopf algebra satisfying the tensor product property.
 Then every thick subcategory of small objects in $K(\Proj\,B)$ is indecomposable.
\end{thm}
 
\begin{proof} 
 The augmentation ideal of $\Ext_B^{*,*}(k,k)$ which consists of all elements in positive degrees is the unique maximal homogeneous ideal and therefore is contained in 
every thick support of $\Spec \Ext_B^{*,*}(k,k)$. Now by theorem \ref{th:grand} it follows that 
every thick subcategory of small objects in this category is indecomposable.
\end{proof}

\begin{cor} Let $B$ denote a finite dimensional sub-Hopf algebra of the mod-2 Steenrod algebra. Then every thick ideal of small objects in  
$K(\Proj\,B)$ is indecomposable.
\end{cor}

\begin{proof} It has been shown (\cite[Corollary 8.6]{hp} and \cite{hp2}) that the tensor product property holds in $K(\Proj\, B)$ when $B$ is a finite dimensional sub-Hopf 
algebra of the mod-2 Steenrod algebra, so corollary follows from the above theorem.
\end{proof}

\section{Stable module categories over Hopf algebras} \label{sec:KS-hopfalgebras}

In this section we give a generalised version of Krause's result \cite[Theorem A]{Kr}. 
To start, let $B$ denote a finite dimensional graded co-commutative Hopf algebra over a field
$k$ of positive characteristic. Following \cite{hp}, we say that $B$ satisfies the tensor product if the category $K(\Proj\,B)$ does. We now prove a 
Krull-Schmidt theorem for $\stmod(B)$, whenever $B$ satisfies the tensor product property.  The category $\StMod(B)$ fails to be a Noetherian stable homotopy category
in general. (In fact, $\pi_*(S)= \Hom_{\stmod(B)}(\Sigma ^* k,k)$ is isomorphic to the Tate cohomology of $B$ which is not Noetherian in general.) So we take the
standard route of using a finite localisation functor \cite[Section 9]{mps} to get around this problem.

We now review some preliminaries from \cite[Section 9]{mps}. There is a finite localisation functor $L_B^f: K(\Proj\,B) \rar K(\Proj\,B)$, whose finite 
acyclics are precisely the objects in the thick subcategory generated by $B$. It has been shown \cite[Theorem 9.6.4]{mps} that the category $\StMod(B)$ is equivalent 
to the full subcategory of $L_B^f$-local objects. This finite localisation functor also establishes a bijection between the nonzero thick ideals of 
small objects in $K(\Proj\,B)$ and the thick ideals of $\stmod(B)$. It follows that the thick ideals of $\stmod(B)$ are in poset bijection with the non-empty specialisation 
closed subsets of $\Spec \, \Ext_{B}^{*,*}(k,k)$. (Note that for the category $K(\Proj\,B)$, the empty set of $\Spec \, \Ext_{B}^{*,*}(k,k)$ corresponds to the thick 
subcategory consisting of the zero chain complex.) Finally let $\mathfrak{m}$ denote the unique bihomogeneous maximal ideal of $\Ext_B^{*,*}(k,k)$ which corresponds to the
trivial thick subcategory in $\stmod(B)$. Then here is the main result of this section.

\begin{thm} \label{thm:mymainthm} Let $B$ denote a finite dimensional graded co-commutative Hopf algebra satisfying the tensor product property. Then every thick ideal 
of $\stmod(B)$ admits a unique Krull-Schmidt decomposition. Conversely, given thick ideals $\T_i$ in $\stmod(B)$ satisfying $\T_i \cap \T_j = 0$ for $i \ne j$,
there exists a thick ideal $\T$ such that $\T \cong \coprod \T_i$. Moreover, 
\[K_0(\T) \cong \underset{i \in I}{\bigoplus}\; K_0(\T_i).\]
\end{thm}

\begin{proof} Recall once again that the thick ideals of the stable module category are in poset bijection with the non-empty subsets of  
$\Spec \; \Ext_{B}^{*,*}(k,k)$ that are closed under specialisation. To begin, let $\T$ be a thick ideal of $\stmod(B)$ and let $S$ denote the 
corresponding specialisation closed subset. Since $\Ext_B^{*,*}(k,k)$ is a graded Noetherian algebra, it follows that 
$S$ admits a unique decomposition $S = \bigcup_i S_i$ such that $S_i \cap S_j = \mathfrak{m}$ 
 for all $i \ne j$ and $S_i$ indecomposable (the proof of
Proposition \ref{prop:graph} can be adapted to the set $S-\mathfrak{m}$ by making use of lemma \ref{le:benson}). 
We now claim that $\T = \coprod \T_{S_i}$, where $\T_{S_i}$ denotes the thick ideal of $\stmod(B)$ corresponding to $S_i$. 

The first thing to note is that for $i \ne j$, $\T_{S_i} \cap \T_{S_j} = 0$:
Suppose not, then the non-zero thick ideal $\T_{S_i} \cap \T_{S_j}$ corresponds to a non-empty subset of $S_i \cap S_j-\mathfrak{m}$ which is impossible.
Therefore it follows that  $\T_{S_i} \cap \T_{S_j} = 0$. 

The next step is to show that the objects in $\T$ are the finite coproducts of the objects in 
$\T_{S_i}$. We give an indirect proof (as before) by showing that $\coprod \T_{S_i}$ is a thick ideal. Then the theorem follows because 
the
specialisation closed subset corresponding to $\coprod \T_{S_i}$ is clearly $\bigcup S_i = S$. To show that $\coprod \T_{S_i}$
is thick (ideal condition holds trivially), by lemma \ref{le:unclear}, all we need to show is that $\Hom_{\stmod}(A,B)= 0$, whenever $A \in \T_{S_i}$
and $B \in \T_{S_j}$ ($i \ne j$).  This was already done for group algebras in proposition \ref{prop:slick}. The proof given there can be easily 
generalised to finite dimensional co-commutative Hopf algebras. So that completes the proof of the existence of Krull-Schmidt decompositions.

For the converse, define $\T$ be the full subcategory of all finite coproducts of objects in $(\T_i)_{i \in I}$ and use the 
fact that $\Hom_{\stmod}(\T_i, \T_j) = 0$ for $i \ne j$, to complete the argument. 

The statement about $K_0$ groups follows from theorem  \ref{th:main}.
\end{proof}

Since tensor product property holds for the finite dimensional sub-Hopf algebras of the mod-2 Steenrod algebra \cite[Corollary 8.6]{hp}, the following
corollary is immediate.

\begin{cor} Let  $B$ denote a finite dimensional sub-Hopf algebra of the mod-2 Steenrod algebra. Then every thick ideal in $\stmod(B)$ admits a Krull-Schmidt
decomposition.
\end{cor}

\begin{rem} It is worth pointing out the difference between the decompositions in theorem \ref{thm:KSforstable(B)}
and those in theorem \ref{thm:mymainthm}. Every thick ideal of small objects in $K(\Proj \, B)$ is indecomposable; on the other hand, thick ideals of small 
objects in $\StMod(B)$ admit Krull-Schmidt decompositions. This striking difference is really the effect of
Bousfield localisation ($\StMod(B)$ is a Bousfield localisation of $K(\Proj \, B)$).
\end{rem}

\subsection{Sub-Hopf algebras of the mod-$2$ Steenrod algebra}

Let $A$ denote the mod-2 Steenrod algebra. We now look at some standard sub-Hopf algebras of $A$ and analyse decompositions in their stable module categories.
For every non-negative integer $n$, denote by $A(n)$ the finite dimensional sub-Hopf algebra of $A$ generated by $\{ \Sq^k |\, 1 \le k \le 2^n\}$.

\vskip 3mm \n
$A(0)$ -- This is the sub-Hopf algebra of $A$ generated by $\Sq^1$ and therefore it is isomorphic to the exterior algebra $\ftwo[\Sq^1]/(\Sq^1)^2$ on $\Sq^1$. It can be easily 
shown that $\Ext_{A(0)}^{*,*} (\ftwo, \ftwo) \cong  \ftwo [h_0]$ with $|h_0|=(1,1)$. Therefore the bihomogeneous prime spectrum of $\Ext_{A(0)}^{*,*} (\ftwo, \ftwo)$
is $\{(0), (h_0)\}$. It is now clear from the thick 
subcategory theorem that there are only two  thick subcategories of $\stmod(A(0))$; the trivial category which corresponds to maximal prime $(h_0)$ and the entire 
category $\stmod(A(0))$ which corresponds to $\{(0), (h_0)\}$.  These subsets of prime ideals are clearly indecomposable closed subsets and therefore we conclude by 
theorem \ref{thm:mymainthm} that all thick subcategories of $\stmod(A(0))$ are indecomposable.
 
\vskip 3mm \n
$A(1)$ -- This is the sub-Hopf algebra of $A$ generated by $\Sq^1$ and $\Sq^2$. This algebra can be shown (using the Adem relations) to be isomorphic to the free 
non-commutative  algebra (over $\ftwo$) generated by symbols $\Sq^1$ and $\Sq^2$ modulo the two sided ideal 
\[(\Sq^1 \Sq^1,\, \Sq^2 \Sq^2 + \Sq^1 \Sq^2 \Sq^1).\] 
The cohomology algebra of $A(1)$ has been computed in \cite{wil}:
\[\Ext_{A(1)}^{*,*}(\mathbb{F}_2,\mathbb{F}_2) \cong \mathbb{F}_2[h_0,h_1,w,v]/(h_0h_1, h_1^3,h_1w,w^2-h_0^2v),\]
where the bi-degrees are: $|h_0|=(1,1), |h_1|=(1,2), |w|=(3,7)$, and $|v|=(4,12)$. The bihomogeneous prime lattice of the above cohomology ring consists of the 
minimal prime $(h_1)$; primes $p_0 = (h_1, w, h_0)$ and $p_2 = (h_1, w, v)$ of height one; and the maximal prime $\mathfrak{m} = (h_1, w, h_0, v)$; see figure \ref{A(1)}.

\begin{figure}
\[
\xymatrix{
                           &  \mathfrak{m} \ar@{-}[dl] \ar@{-}[dr] &  \\
p_0  \ar@{-}[dr]  &                                                &  p_1 \ar@{-}[dl] \\
                           &  (h_1)                                &                                  
}
\]
\label{A(1)}
\caption{Bihomogeneous prime lattice of $\Ext_{A(1)}^{*,*} (\ftwo, \, \ftwo)$}
\end{figure}
The thick subcategory theorem now tells us that
there are five thick subcategories in $\stmod(A(1))$ which correspond to the thick supports $\{\mathfrak{m}\}$, $\{p_0, \mathfrak{m}\}$, 
$\{p_1, \mathfrak{m}\}$, $\{p_0, p_1, \mathfrak{m}\}$, and $\{(h_1), p_0, p_1, \mathfrak{m}\}$. It is clear that the only decomposable thick support in this 
list is $\{p_0, p_1, \mathfrak{m}\}$; in fact,
 \[\{p_0, p_1, \mathfrak{m} \} = \{ p_0, \mathfrak{m}\} \, \cup \,\{ p_1, \mathfrak{m}\}\] 
is the desired decomposition. This induces a Krull-Schmidt decomposition of thick subcategories supported on these thick supports:
 \[ \T_{\{p_0, p_1, \mathfrak{m}\}} \cong  \T_{\{p_0, \mathfrak{m}\}} \amalg \T_{\{p_1, \mathfrak{m}\}}.\] 
Note that the remaining four thick subcategories in $\stmod(A(1))$ are indecomposable.

\vskip 3mm 
The analysis of decompositions in stable module categories over $A(2)$ and higher are much more complicated.

\section{Finite spectra} \label{sec:KS-finitespectra}

In this section we work in the triangulated category $\F$ of finite spectra. Some good references for spectra are \cite{bluebook} and \cite{mar}.
Recall that a  spectrum $X$ is 
torsion  if $\pi_*(X)$  is a torsion group. In particular, a finite spectrum $X$ is torsion ($p$-torsion) if $\pi_*(X)$ consists of finite abelian groups 
($p$-groups) in every degree. We will now show that the thick subcategories of  finite spectra admit Krull-Schmidt decompositions. 

$\FT$ will denote the category of finite torsion spectra, and 
$\C_{n,p}$ will denote the thick subcategory of finite $p$-torsion spectra consisting of all $K(n-1)_*$ acyclics, where $K(n-1)$ denotes the $(n-1)$st Morava 
K-theory at the prime $p$.

\begin{lemma} If $\A$ is a thick subcategory of $\F$, then either $\A = \F$ or $\A \subseteq \FT$.
\end{lemma}

\begin{proof} It suffices to show that if $X$ is any non-torsion finite spectrum, then the integral sphere spectrum $S$ belongs to the thick 
subcategory generated by $X$. Consider the cofibre sequence
\[ X \stk{p} X \rar X \wedge M(p) \rar \Sigma X. \]
It is easy to see that $X \wedge M(p)$ is a type-1 spectrum, therefore the thick subcategory generated by $X$ contains all finite torsion spectra.
The next step is to show that there is a cofibre sequence
\[ \coprod S^r \rar X \rar T \rar \coprod S^{r+1},\]
where $T$ is a finite torsion spectrum; see \cite[Proposition 8.9]{mar}. It is clear from these exact triangles that $S$ belongs to the thick subcategory generated by 
$X$.
\end{proof}

We now start analysing decompositions for thick subcategories of finite spectra.

\begin{prop} The category $\F$ is indecomposable.
\end{prop}
\begin{proof} If $\F = \A \amalg \B$ is any decomposition of $\F$, then that would giving a splitting $A \vee B$ of the sphere spectrum 
($S$).
Since $S$ is an indecomposable spectrum, either $A$ or $B$ must be the sphere spectrum, which implies that either $\A$ or $\B$ is equal to $\F$. 
\end{proof}

We now show that the thick subcategories of finite \emph{torsion} spectra admit Krull-Schmidt decompositions.

\begin{thm}  Every thick subcategory $\D$ of $\FT$ admits a unique Krull-Schmidt decomposition.
Conversely, given any collection of thick subcategories $\D_i$ in $\FT$ such that $\D_i \cap \D_j = 0$ for all $i \ne j$, 
there exits a  thick subcategory $\D$ such that $\underset{i}{\coprod}\D_i$ is a Krull-Schmidt decomposition for $\D$.
Moreover $K_0(\D) \cong \bigoplus_i K_0(\D_i).$
\end{thm}

\begin{proof}  
For every prime number $p$, define $\D_p := \D \cap \C_{1,p}$. By \cite[8.2, Theorem 20]{mar}  every torsion spectrum $X$ can be written as   
$X = \bigvee_p X_{(p)}$ where $X_{(p)}$ is the $p$-localisation of $X$. If $X$ is a finite torsion spectrum, it is easy to see that $X_{(p)}$ 
is a finite $p$-torsion spectrum and therefore belongs to $\C_{1,p} (\subseteq \FT)$. Since $\D$ is a thick subcategory of $\FT$, all the 
summands $X_{(p)}$ also belong to $\D$. Therefore $X_{(p)}$ belongs to $ \C_{1,p} \cap \D =\D_p$.
Also since $X$ is a finite spectrum, the above wedge runs over only a finite set of primes. This implies that the spectra in 
$\D$ are the finite coproducts of spectra in $(\D_p)$.

$\D_p$ is a thick subcategory in $\FT$:  Showing the $\C_{1,p}$ is a thick subcategory of $\FT$ is a routine verification. Since intersection of
thick subcategories is thick,  $\D \cap \C_{1,p}$ is thick in $\FT$.

$\D_p \cap \D_q = 0:$   It is clear that if a finite spectrum is both $p$-torsion and $q$-torsion, for primes $p \ne q$, then  
it  is trivial, i.e., $\C_{1,p} \cap \C_{1,q} = 0$. But $\D_p \cap \D_q \subseteq \C_{1,p} \cap \C_{1,q}$, therefore  $\D_p \cap \D_q = 0.$

$\D_p$ is indecomposable:  $\D_p$ is by definition a thick subcategory of $\C_{1,p}$.  The thick subcategories of $\C_{1,p}$ are well-known to be
nested \cite{hs}. In particular they are all indecomposable.

This completes the proof of the existence of a Krull-Schmidt decomposition for any thick subcategory $\D$ in $\FT$. 

Uniqueness: Let $\coprod \E_k$ be another Krull-Schmidt decomposition for $\D$. Then since each $\E_k$ is an indecomposable thick
subcategory of $\FT$, we conclude that $\E_k = \C_{n_k, p_k}$ for some integer $n_k$ and some prime $p_k$. Thus the decomposition is unique.

For the converse, define $\D$ be the full subcategory of all finite coproducts of objects from $\D_p$. We just have to show that $\D$
is a thick subcategory and then the rest follows.  Note that $[X,Y] = 0$ whenever $X$ is a $p$-torsion spectrum and $Y$ is a $q$-torsion spectrum \cite[8.2, Lemma 21]{mar}. 
Therefore thickness of $\D$ follows from lemma \ref{le:unclear}.
\end{proof}

The following corollary is clear from these decompositions.

\begin{cor} Let $\A$ be be a thick subcategory of $\F$. Then either $\A = \F$, or $\A = \coprod_{p \in S} \C_{i_p,p}$, where
S is a subset of prime numbers, and $i_p$ is a positive integer for each prime $p$. Further, in the latter case, we have
\[ K_0(\A) \cong \underset{p \, \in \,S}{\bigoplus} \; K_0(\C_{i_p, p}).\]
\end{cor}

\section{Questions}

\subsection{Krull-Schmidt decompositions in $D^b(\proj\,R)$}

We have shown that the thick subcategories of perfect complexes over a \emph{Noetherian} ring admit Krull-Schmidt
decompositions. Note that proposition \ref{prop:graph} was crucial for  this decomposition. Now the question that 
remains to be answered  is whether  this proposition holds for arbitrary commutative rings. In other words,
if $R$ is any commutative ring, is it true that every thick support of $\Spec(R)$ admits a unique Krull-Schmidt 
decomposition? An affirmative answer to this question would imply, by theorem \ref{thm:aux}, that thick
subcategories of perfect complexes over commutative rings (not necessarily Noetherian) admit Krull-Schmidt decompositions.

\subsection{The tensor product property}

Recall that our decompositions for the thick subcategories of small objects in $K(\Proj\,B)$ and $\StMod(B)$ assumed that 
the finite dimensional co-commutative Hopf algebra $B$ satisfies the tensor product property.  Now the question
that arises is which finite dimensional Hopf algebras satisfy the tensor product property? This seems to be a question of independent interest; 
see \cite[Corollary 3.7]{hp} for some interesting consequences of the tensor product property.
Friedlander and Pevtsova \cite{FriPev} have shown recently that the finite dimensional co-commutative 
\emph{ungraded} Hopf algebras satisfy the tensor product property. Using their ideas, or otherwise, we would like to know if 
this property also holds for finite dimensional graded co-commutative Hopf algebras. Also we do not know if the tensor product product 
property in \cite{FriPev} is equivalent to the one given in \cite{hp}.

\bibliographystyle{alpha}

\end{document}